\documentclass[opre,nonblindrev]{informs3} 

\DoubleSpacedXI 

\usepackage{endnotes}
\let\footnote=\endnote
\let\enotesize=\normalsize
\def\notesname{Endnotes}%
\def\makeenmark{$^{\theenmark}$}
\def\enoteformat{\rightskip0pt\leftskip0pt\parindent=1.75em
  \leavevmode\llap{\theenmark.\enskip}}

\usepackage{multirow}
\usepackage{caption}
\usepackage{subcaption}
\usepackage{epsf}
\usepackage{slashbox}

\usepackage{natbib}
 \bibpunct[, ]{(}{)}{,}{a}{}{,}%
 \def\bibfont{\small}%
 \def\bibsep{\smallskipamount}%
 \def\bibhang{24pt}%
\TheoremsNumberedThrough     
\ECRepeatTheorems

\EquationsNumberedThrough    

\begin{document}


\RUNAUTHOR{Qian, Wang and Wen}

\RUNTITLE{A Composite Risk Measure Framework for Decision Making under Uncertainty}

\TITLE{A Composite Risk Measure Framework for Decision Making under Uncertainty}

\ARTICLEAUTHORS{%
\AUTHOR{Pengyu Qian}
\AFF{School of Mathematical Sciences, Peking University, China, \EMAIL{pengyu.qian@pku.edu.cn}} 
\AUTHOR{Zizhuo Wang}
\AFF{Department of Industrial and Systems Engineering, University of Minnesota, Minneapolis, MN 55455, \EMAIL{zwang@umn.edu}}
\AUTHOR{Zaiwen Wen}
\AFF{Beijing International Center for Mathematical Research, Peking University, China, \EMAIL{wenzw@math.pku.edu.cn}}
} 

\ABSTRACT{%
In this paper, we present a unified framework for decision making
under uncertainty. Our framework is based on the composite of two
risk measures, where the inner risk measure accounts for the risk of
decision given the exact distribution of uncertain model parameters,
and the outer risk measure quantifies the risk that occurs when
estimating the parameters of distribution. We show that the model is
tractable under mild conditions. The framework is a generalization
of several existing models, including stochastic programming,
robust optimization, distributionally robust optimization, etc.
Using this framework, we study a few new models which imply probabilistic
guarantees for solutions and yield less conservative results comparing
to traditional models. Numerical experiments are
performed on portfolio selection problems to demonstrate the
strength of our models.

}%



\maketitle

%


\section{Introduction}\label{sec:intro}
In this paper, we consider a decision maker who wants to minimize an objective function
$H(\BFx, \BFxi)$,
where $\BFx \in \mathbb{R}^n$ is the decision variable and $\BFxi \in \mathbb{R}^s$ is some uncertain/unknown parameter related to the model.
For example, in a newsvendor problem, $\BFx$ is the order amount of newspapers by a newsvendor, and $\BFxi$ is the uncertain future demand. Similarly, in an investment problem, $\BFx$ is the portfolio chosen by a portfolio manager, and $\BFxi$ is the unknown future returns of the instruments.
The existence of the uncertain parameters distinguishes the problem from ordinary optimization problems and has led to several decision making paradigms.

One of the earliest attempts to deal with such decision making
problems under uncertainty was proposed by
\cite{Danzig_linear_uncertainty}, where it was assumed that the
distribution of $\BFxi$ is known exactly and the decision is chosen
to minimize the expectation of $H(\BFx, \BFxi)$. Such an approach is
called stochastic programming. Another approach named robust
optimization initiated by \cite{RO_origin_soyster} supposes that all
possible values of $\BFxi$ lie within an uncertainty set, and the
decision should be made to minimize the worst-case value of $H(\BFx,
\BFxi)$. Stochastic programming and robust optimization models can
be viewed as two extremes in the spectrum of available information
in decision making under uncertainty. There are models in between
these two extremes. For example, distributionally robust models
(\citealt{DRSP_scarf_1958}, \citealt{DRSP_dupacova_1987},
\citealt{Delage-Ye-2010}, etc.) take into account both the
stochastic and the robust aspects, where the distribution of $\BFxi$
is assumed to belong to a certain distribution set and the
worst-case expectation of $H(\BFx, \BFxi)$ is minimized. There are
also various models which minimize certain risk
(\citealt{Rockafellar-2000}, \citealt{gaivoronski2005value}, etc.)
or the worst-case risk (\citealt{worst_case_VaR},
\citealt{Zhu-Fukushima-2009}, etc.) of $H(\BFx, \BFxi)$. We will
present a more detailed review of these models in Section \ref{sec:review}.
In addition to the study of individual models, there have been recent efforts to seek connections between different models and to put forward more general models. For instance, \cite{Bertsimas-uncertainty-2009} and \cite{natarajan2009constructing} show that uncertainty sets can be constructed according to decision maker's risk preference, \cite{Bertsimas-data-driven-2014} propose a general framework for data-driven robust optimization, and \cite{wiesemann2013distributionally} propose a framework for distributionally robust models.

While some models have demonstrated their effectiveness in practice, there are still some ignored issues in the existing literature:
\begin{enumerate}
\item There lacks a unified framework which includes all the models above, namely, stochastic programming model, robust optimization model, distributionally robust model and worst-case risk models.

\item Though risk measures have been imposed on the objective function to deal with parameter uncertainty (\citealt{coherent_risk_measure}, \citealt{Bertsimas-uncertainty-2009}, etc.), no attempt has been made to impose risk measure on the expectation or other functionals of the objective function with respect to distributional uncertainty.

\item Bayesian approach has not been fully considered when modeling decision making under distributional uncertainty in the existing literature, albeit its appropriateness for such problems.
\end{enumerate}

The goal of our paper is to fill these gaps and build a unified modeling framework for decision making under uncertainty. Our unified framework is based on a risk measure interpretation for robustness and encompasses several popular models, such as stochastic programming, robust optimization, distributionally robust optimization, worst-case risk models, etc.
Specifically, we minimize the composite of two risk measures, where the inner risk measure accounts for the risk of decision given the exact distribution of parameters, and the outer risk measure quantifies the risk that occurs when estimating the parameters of distribution.
For the outer risk measure, we take a Bayesian approach and consider the posterior distribution of distribution parameters. We show that the composite of risk measures is convex as long as both the inner and the outer risk measures are convex risk measures. We also use this framework to construct several new models which have real world meanings and perform numerical tests on these models which demonstrate their strength.

We summarize our contributions as follows:
\begin{enumerate}
\item We propose a composite risk measure (CRM) framework for decision making under uncertainty where the composite of two risk measures is minimized. It is a generalization of several existing models. We show that the corresponding optimization problem is convex under mild conditions.
\item We take a novel approach to deal with distributional uncertainty by making use of risk measure and Bayesian posterior distribution of distribution parameters.
\item Using the composite risk measure framework, we study a VaR-Expectation model, a CVaR-Expectation model and a CVaR-CVaR model, investigating their tractability and probabilistic guarantees. Numerical experiments show that these models can be solved efficiently and perform well in portfolio selection problems.
\end{enumerate}

The remainder of this paper is organized as follows. In Section \ref{sec:review}, we briefly review the existing models for decision making under uncertainty. In Section \ref{sec:CRM_frame}, we present our composite risk measure framework and show that several existing models fall into the framework. Several new models within the general framework are proposed in Section \ref{sec:new_model}. Numerical results for these new models are shown in Section \ref{sec:numer}.

\textbf{Notations.} Throughout the paper, the following notations
will be used. Ordinary lower case letters $(x,y,\dots)$ denote scalars,
boldfaced lower case letters $(\BFx,\BFy,\dots)$ denote vectors. Specifically,
$\BFx \in \mathbb{R}^n$ denotes the decision variable and
$\BFxi \in \mathbb{R}^s$ denotes the uncertain/unknown parameters.
$H(\BFx,\BFxi)$ is the loss function under decision $\BFx$ and parameter $\BFxi$.

\section{Review of Existing Models for Decision Making under Uncertainty}\label{sec:review}
There have been many models proposed in the literature that study decision making problem under uncertainty.
In this section, we provide a review of those existing models.
\subsection{Stochastic Programming}\label{subsec:review_SP_model}
One of the most popular ways to solve decision making problem under uncertainty is through stochastic programming. In stochastic programming models, one assumes that the full distributional information of $\BFxi$ is available. Then to choose an optimal decision variable $\BFx$, one considers a certain functional of the random loss $H(\BFx, \BFxi)$. Examples of this approach include:

\begin{descr}
\item[\textbf{Expectation optimization.}] \cite{Danzig_linear_uncertainty} considers the case where the objective is to minimize the expectation of the random loss $H(\BFx, \BFxi)$. Namely, the optimization problem is
    \begin{eqnarray}\min_{\BFx}&&\mathbb{E} [H(\BFx, \BFxi)].\label{prob_opt_expectation}\end{eqnarray}

    There are much literature that study such optimization problems. We refer interested readers to \cite{shapiro2009lectures} for a comprehensive review.

\item[\textbf{Value-at-Risk (VaR) optimization.}] VaR has enjoyed great popularity as a measure of downside risk in finance industry since its introduction in the 1990s (see \citealt{JPMorgan_VaR_1996}).
    For fixed $\BFx$, the $\delta$-VaR of $H(\BFx, \BFxi)$ is defined as the $\delta$ quantile of the random loss $H(\BFx, \BFxi)$. Mathematically,
    \begin{eqnarray} \text{VaR}_{\delta}(H(\BFx, \BFxi)) & \triangleq & \inf \left\{t \in \mathbb{R}|\mathbb{P}(H(\BFx, \BFxi) \geqslant t) \leqslant 1-\delta \right\}\label{VaR_definition}.\end{eqnarray}
    The corresponding VaR optimization problem can be written as:
    \begin{eqnarray}
    \min_{\BFx, t}&&t  \label{prob_opt_VaR}\\
    \mbox{s.t.}&&\mathbb{P}(H(\BFx, \BFxi) \geqslant t) \leqslant 1-\delta. \nonumber
    \end{eqnarray}

    The VaR optimization problem has been studied extensively in the literature, see, e.g., \cite{kast1998var}, \cite{lucas1998extreme}, etc.
    However, there are three main drawbacks of VaR. First, it doesn't take into account the magnitude of loss beyond VaR, resulting in decision maker's preference for taking `excessive but remote' risks (\citealt{VaR_drawback_ignoretail}). Second, it is not subadditive (\citealt{coherent_risk_measure}), meaning that the VaR of a combined portfolio can be larger than the sum of the VaRs of its components, a property that is not desired in many applications.
    Lastly,
    $\text{VaR}_{\delta}(H(\BFx, \BFxi))$ is usually not convex in $\BFx$ (\citealt{birge_stochastic_intro_2002}), making the optimization problem intractable.

\item[\textbf{Conditional Value-at-Risk (CVaR) optimization.}] To overcome the drawbacks of VaR, researchers further proposed a modified version of VaR, the CVaR (also called the expected shortfall in some literature).
    For a certain random loss $X$, the $\delta$-CVaR is the expected loss in the worst $1-\delta$ cases. Mathematically, it can be written as:
    \begin{eqnarray*}
    \textup{CVaR}_{\delta}(X) &=& \frac{1}{1-\delta} \int_{\delta}^1 \textup{VaR}_s(X)ds.
    \end{eqnarray*}
    For atomless distributions, CVaR is equivalent to the conditional expectation of the loss beyond VaR, namely,
    \begin{eqnarray*}\textup{CVaR}_{\delta}(X) &=& \mathbb{E}[X|X\geqslant \textup{VaR}_{\delta}(X)].\end{eqnarray*}
    \cite{Rockafellar-2000} show that CVaR can be obtained by solving the following convex program:
    \begin{equation}
    \min_{\alpha \in \mathbb{R}}\ \alpha + \frac{1}{1-\delta} \mathbb{E}[(H(\BFx, \BFxi)- \alpha)^{+}],\label{prob_opt_CVaR}
    \end{equation}
    which leads to another definition of CVaR. This formulation brings computational convenience as the objective function is explicit and convex in $\alpha$.

\end{descr}
\subsection{Robust Optimization}\label{subsec:review_RO_model}
Another popular way for decision making under uncertainty is to use a robust optimization approach.
In the robust optimization approach, instead of assuming the distributional knowledge of $\BFxi$, one only assumes that $\BFxi$ takes values in a certain uncertainty set $\Xi$. Then, when choosing the decision variable, one considers the worst-case outcome associated with each decision, where the worst-case scenario is chosen from the specified uncertainty set. The optimization problem can be written as follows:
\begin{equation}
\min_{\BFx} \ \max_{\BFxi \in \Xi} \ H(\BFx, \BFxi). \label{prob_opt_plain_robust}
\end{equation}

Robust optimization problems have more general forms where constraints instead of objective function are affected by parameter uncertainty (e.g. \citealt{Bertsimas_review_RO}).
Robust optimization was first proposed by \cite{RO_origin_soyster} and has attracted much attention in the past few decades. For a comprehensive review of the literature, we refer readers to \cite{Bertsimas_review_RO} and \cite{RO_book_BenTal}.

Choosing a suitable uncertainty set is essential in formulating a robust optimization problem. Two main issues should be taken into consideration when designing uncertainty sets:

\begin{descr}
\item[\textbf{Tractability.}] Only a few uncertainty sets will lead to tractable counterparts for the original problem. Some known cases include polyhedral uncertainty set (\citealt{ellipsoidal_RO_BenTal}), ellipsoidal uncertainty set (\citealt{ellipsoidal_RO_BenTal}, \citealt{ellipsoidal_RO_Ghaoui_1} and \citealt{ellipsoidal_RO_Ghaoui_2}), norm uncertainty set (\citealt{norm_RO_Bertsimas}), etc.
\item[\textbf{Conservativeness.}] Intuitively, when the uncertainty set is very large, the resulting decision will be very robust, but sometimes too conservative to be of practical use. Much work has been done on the choice of uncertainty set, see \cite{price_RO_Bertsimas}, \cite{Chen-2007} and \cite{Bertsimas-data-driven-2014}.
\end{descr}

\subsection{Distributionally Robust Models}\label{subsec:review_DRO_model}
In practice, it is often the case that one only has partial information about the distribution of $\BFxi$, such as first and second moments. Applying stochastic programming in those cases is not feasible. Meanwhile, using a robust optimization approach is not straightforward either and often results in overly-conservative solutions.
As a result, an intermediate path has been proposed which is called the distributionally robust optimization (DRO) model. In the DRO models, the decision maker constructs an uncertainty set $\mathcal{F}$ of the underling distribution and minimizes the expected value of the objective function under the worst-case distribution chosen from $\mathcal{F}$. That is, one considers the following problem:
\begin{equation} \min_{\BFx}\ \max_{F \in \mathcal{F}}\ \mathbb{E}_{\BFxi\sim F} [H(\BFx, \BFxi)], \label{distributionally_robust}\end{equation}
where $\mathbb{E}_{\BFxi\sim F} [\cdot]$ denotes the expectation of $H(\BFx, \BFxi)$ when $F$ is the distribution of $\BFxi$.

Distributionally robust model was first proposed in \cite{DRSP_scarf_1958}. After its introduction, several different choices of the uncertainty set $\mathcal{F}$ have been proposed. We will review two most discussed types of uncertainty sets below and refer readers to \cite{wiesemann2013distributionally} and references therein for other types of uncertainty sets.

\begin{descr}
\item[\textbf{Distributions with partial knowledge on moments.}]\cite{DRSP_scarf_1958}, \cite{DRSP_dupacova_1987}, \cite{prekopa1995stochastic}, \cite{optimal_inequality_probability}, etc. consider a family of distributions with known moments.
    Specifically, \cite{DRSP_scarf_1958} consider the following optimization problem (in the context of an inventory problem):
    \begin{align}
    &\min_{x \in \mathbb{R}} \ \max_{F \in \mathcal{F}{(\mu_0,\sigma_0})}\
    \mathbb{E}_{\xi\sim F} [c(x - \xi)^{+} + r(\xi - x)^{+}] \nonumber\\
    &\text{where }
    \mathcal{F}(\mu_0, \sigma_0) = \left\{
    F(\xi)\in \mathcal{M} \left|
    \begin{aligned}
    &\mathbb{P}(\xi \in \mathbb{R}^{+}) = 1 \\
    &\mathbb{E}_{\xi\sim F}(\xi) = \mu_0 \\
    &\mathbb{E}_{\xi\sim F}(\xi - \mu_0)^2 = \sigma_0^2
    \end{aligned}
    \right\}\right. \label{prob_scarf} .
    \end{align}
    In (\ref{prob_scarf}), $\mathcal{M}$ is the set of all probability measures in the probability space where $\xi$ is defined and $c$ and $r$ are given parameters. \cite{DRSP_scarf_1958} show that the worst-case distribution is a two-point distribution with given decision $x$ and derive a closed-form solution to (\ref{prob_scarf}).

    \cite{Delage-Ye-2010} consider a more general form of moment constraints. Namely, the uncertainty set is constructed as:
    \begin{equation}
    \mathcal{F} (\Xi,\BFmu_0,\Sigma_0,\gamma_1,\gamma_2) = \left\{
    F(\BFxi)\in \mathcal{M}  \left|
    \begin{aligned}
    &\mathbb{P} (\BFxi \in \Xi)=1 \\
    &(\mathbb{E}_{\BFxi\sim F}(\BFxi)-\BFmu_0)^T \Sigma_0^{-1} (\mathbb{E}_{\BFxi\sim F}(\BFxi)-\BFmu_0)\leqslant \gamma_1 \\
    &\mathbb{E}_{\BFxi\sim F} [(\BFxi-\BFmu_0)(\BFxi-\BFmu_0)^T] \preceq \gamma_2 \Sigma_0
    \end{aligned}
    \right\}.\right. \label{delage_ye}
    \end{equation}
    \cite{Delage-Ye-2010} show that problem (\ref{delage_ye}) can be solved by a convex optimization problem. They also provide a data-driven method for choosing the parameters $\gamma_1$ and $\gamma_2$.

    More recently, \cite{wiesemann2013distributionally} consider a general model which incorporates both moment information and support information. The distributional uncertainty set considered is:
    \begin{equation}
     \mathcal{F} = \left\{
     F \in \mathcal{M}(\mathbb{R}^l \times \mathbb{R}^m) \left|
     \begin{aligned}
     &\mathbb{E}_{\BFxi\sim F}(A \BFxi + B \BFeta) = b\\
     &\mathbb{P}_{\BFxi\sim F}[(\BFxi,\BFeta)\in \mathcal{C}_{i}]\in [\bar{\BFp}_{i},\underline{\BFp}_{i}],\forall i\in\mathcal{I}
     \end{aligned}
     \right\}\right. \label{wiesemann_framework},
    \end{equation}
    where $\mathcal{M}(\mathbb{R}^l \times \mathbb{R}^m)$ represents probability distribution on $\mathbb{R}^l \times \mathbb{R}^m$, $\BFxi \in \mathbb{R}^{l}$ is the random term, $\BFeta \in \mathbb{R}^{m}$ is an auxiliary random vector, $A\in\mathbb{R}^{k\times l}$, $B\in\mathbb{R}^{k\times m}$, $b\in\mathbb{R}^{k}$ are predetermined parameters and $\mathcal{C}_{i}$ are conic confidence sets. It is shown that when certain conditions are satisfied, model (\ref{distributionally_robust}) with distributional uncertainty set (\ref{wiesemann_framework}) can be solved by conic programming.
\item[\textbf{Distributions related to a given distribution.}] \cite{Ben-Tal-2013} consider a set of distributions that arise from $\phi$-divergence. Namely,
    \begin{equation}
    \mathcal{F} = \left\{
    f \in \mathbb{R}^m \left|I_{\phi}(f,\hat{f}) \leqslant \rho, \sum_{i=1}^m f_i = 1,f_i \geqslant 0, i=1,2,\dots,m\right\}\right. \label{bental_divergence},
    \end{equation}
    where $\hat{f}$ is a given probability vector (e.g., the empirical distribution vector). $I_{\phi}(f,\hat{f})$ is the $\phi$-divergence between two probability vectors $f$ and $\hat{f}$ (\citealt{pardo2005statistical}). \cite{Ben-Tal-2013} show that the robust linear constraint
    $$(\BFa + \mathbb{E}_{\BFxi\sim F}(\BFxi))^{T}\BFx \leqslant \beta, \quad \forall F \in \mathcal{F},$$
    can be written equivalently as a linear constraint, a conic constraint, or a convex constraint,
    depending on the choice of $\phi$, where $\BFa \in \mathbb{R}^n$ and $\beta \in \mathbb{R}$ are fixed parameters,
    $\BFx \in \mathbb{R}^n$ is the decision variable and $\mathcal{F}$ is defined in (\ref{bental_divergence}). Similar approaches have been taken by \cite{Z.Wang-2013} and \cite{klabjan2013robust_chi}.

    \cite{Bertsimas-data-driven-2014} propose a model where the distributional uncertainty set is constructed by means of hypothesis test given a set of available data. Namely, two hypotheses are compared:
    $$H_0:\BFp^{*}=\BFp_0\quad vs.\quad H_A:\BFp^{*}\neq\BFp_0,$$
    where $H_0$ is the null hypothesis, $H_A$ is the alternative hypothesis and $\BFp_0$ is an arbitrary distribution. By specifying a hypothesis test, e.g. $\chi^2$-test, G-test, etc., and a confidence level $\epsilon$, one can construct a distribution set $\mathcal{F}$ containing all the distributions that pass the test under the given set of data.
\end{descr}

\subsection{Other Choices of Objectives}\label{subsec:review_other_model}
In addition to the models mentioned above, there are several other models for decision making under uncertainty that have been studied in the literature.
\begin{descr}
\item[\textbf{Minimizing worst-case VaR.}] \cite{worst_case_VaR} consider the problem of minimizing VaR (defined in (\ref{VaR_definition})) over a portfolio of random loss $\BFxi$, where only partial knowledge about the distribution $F$ of $\BFxi$ is known. Mathematically, the optimization problem is:
    \begin{eqnarray}
    \min_{\BFx,t}&&t  \nonumber\\
    \mbox{s.t.}&&\mathbb{P}_{\BFxi\sim F}(\BFxi^{T}\BFx \geqslant t) \leqslant 1-\delta,
    \quad \forall F \in \mathcal{F} \label{worst_case_VaR} \\
    &&\BFx \in \mathcal{X}, \nonumber
    \end{eqnarray}
    where $\mathcal{F}$ is the set of all probability measures with given first two moments $\mu_0$ and $\Sigma_0$ ($\Sigma_0 \succeq 0$). Using the exact Chebyshev bound given in \cite{optimal_inequality_probability}, \cite{worst_case_VaR} show that problem (\ref{worst_case_VaR}) can be reduced to a second-order cone program (SOCP). \cite{worst_case_VaR} further show that if $\mathcal{F}$ is a set of distributions with the first and second moments $(\mu,\Sigma)$ satisfying
    $(\mu,\Sigma) \in \text{conv}\left\{(\mu_1,\Sigma_1),(\mu_2,\Sigma_2),\dots,(\mu_k,\Sigma_k)\right\}$
    or
    $(\mu,\Sigma)$ is bounded in a componentwise fashion,
    problem (\ref{worst_case_VaR}) can be solved by an SOCP or a semi-definite program (SDP) respectively.

\item[\textbf{Minimizing worst-case CVaR}]\cite{Zhu-Fukushima-2009} solve the problem of optimizing the CVaR (defined in (\ref{prob_opt_CVaR})) of a portfolio when the distribution $F$ of the random return $\BFxi$ is only known to belong to an uncertainty set $\mathcal{F}$ instead of being exactly known, namely,
    \begin{equation} \min_{\BFx \in \mathcal{X}}\ \max_{F \in \mathcal{F}}\ \textup{CVaR}_{F,\delta} [H(\BFx,\BFxi)],
    \label{worst_case_CVaR}\end{equation}
    where $\textup{CVaR}_{F,\delta}(\cdot)$ denote the $\delta$-CVaR of a random variable whose distribution is $F$.
    \cite{Zhu-Fukushima-2009} show that for certain forms of $\mathcal{F}$, problem (\ref{worst_case_CVaR}) can be transformed to a convex optimization problem.
\end{descr}

\section{A Composite Risk Measure Framework}\label{sec:CRM_frame}
In this section, we present a unified framework for decision making problem under uncertainty. Our framework encompasses all the decision paradigms discussed in Section \ref{sec:review} and can be used to generate new ones.

The idea of our framework is based on risk measures defined as follows:
\begin{definition}
Let $\mathcal{L}$ be a set of random variables defined on the sample space $\Omega$. A functional $\rho(\cdot): \mathcal{L} \rightarrow \mathbb{R}$ is a risk measure if it satisfies the following properties:
\begin{enumerate}
\item \emph{Monotonicity:} For any $\BFX, \BFY \in \mathcal{L}$, if $\BFX \geqslant \BFY$, then $\rho(\BFX) \geqslant \rho(\BFY)$, where $\BFX\geqslant \BFY$ means that $\BFX(\omega) \geqslant \BFY(\omega)$ for any $\omega \in \Omega$.
\item \emph{Translation invariance:} For any $\BFX \in \mathcal{L}$ and $c \in \mathbb{R}$, $\rho(\BFX+c)=\rho(\BFX)+c$.
\end{enumerate}
\end{definition}

In addition, \cite{coherent_risk_measure} define a subset of risk measures satisfying some structural properties presented as follows.
\begin{definition}
If a risk measure $\rho(\cdot)$ satisfies the following properties:
\begin{enumerate}
\item \emph{Convexity:} For any $\BFX, \BFY\in \mathcal{L}$ and $\lambda \in [0,1]$, $\rho(\lambda\BFX+(1-\lambda)\BFY) \leqslant \lambda \rho(\BFX)+(1-\lambda)\rho(\BFY)$;
\item \emph{Positive homogeneity:} For any $\BFX \in \mathcal{L}$ and $\lambda\geqslant 0$, $\rho(\lambda \BFX)=\lambda\rho(\BFX)$,
\end{enumerate}
then $\rho(\cdot)$ is called a coherent risk measure. If only convexity holds, it is called a convex risk measure (\citealt{follmer2002convex}).
\end{definition}

To establish our framework, first we note that given $\BFx$, $H(\BFx, \BFxi)$ is a random variable defined on the sample space of $\BFxi$ (which we denote by $\Omega_0$).  We denote such a random variable by $Y(\BFx)$.
Define $g_F(\cdot)$ to be a risk measure for $Y(\BFx)$, where the subscript $F$ is used to show the dependence of this risk measure on the choice of distribution $F$.
Now we further define a measurable space for $F$: $(\Omega_1, \Sigma_1)$, where $\Omega_1$ denotes the space of all the distribution functions for $\BFxi$, and $\Sigma_1$ is a $\sigma$-algebra defined on the space of such distributions. Moreover, we can define a measure $\mathbb{P}_1$ on such a space using concepts from Bayesian statistics (see the following passage
for detailed discussions). With this definition, the risk measure $g_F(Y(\BFx))$ can be viewed as a random variable too in the following way:
\begin{eqnarray}
Z(\BFx): F\in\Omega_1 &\rightarrow& g_F(Y(\BFx)) \in {\mathbb R}.
\end{eqnarray}
We denote the linear space of $Z(\BFx)$ by $\mathcal{Z}$.
Finally, since $Z(\BFx)$ is a random variable, we can apply another risk measure $\mu: \mathcal{Z} \rightarrow \mathbb{R}$ and consider the following optimization problem:
\begin{equation}\min_{\BFx \in \mathcal{X}}\ \mu(Z(\BFx)). \label{composite_naive} \end{equation}
Therefore, our framework can be written as follows:
\begin{equation}
\min_{\BFx \in \mathcal{X}}\ \mu(g_{F}(H(\BFx,\BFxi))). \label{coherent_problem_1}
\end{equation}
We call this the composite-risk-measure (CRM) framework for decision making under uncertainty.
In the following discussions, we will refer to $g_{F}(\cdot)$ as the inner risk measure and $\mu(\cdot)$ as the outer risk measure. We first present the following tractability result for problem (\ref{coherent_problem_1}).

\begin{proposition}\label{coherent_composite_convex}
Optimization problem (\ref{coherent_problem_1}) is a convex optimization problem if the following holds:
\begin{enumerate}
\item $H(\BFx,\BFxi)$ is convex in $\BFx$;
\item $\mathcal{X}$ is convex;
\item $\mu(\cdot)$ is a convex risk measure;
\item $g_{F}(\cdot)$ is a convex risk measure for each $F\in\mathcal{F}$.
\end{enumerate}
\end{proposition}

\proof{Proof.}
We show that $\mu(g_{F}(H(\BFx,\BFxi)))$ is a convex function of $\BFx$. For any $\BFx, \BFy \in \mathcal{X}$ and $0\leqslant \lambda \leqslant 1$, we have:
\begin{eqnarray*}
\lambda \mu(g_{F}(H(\BFx,\BFxi))) + (1-\lambda) \mu(g_{F}(H(\BFy,\BFxi)))
&\geqslant& \mu(\lambda g_{F}(H(\BFx,\BFxi)) + (1-\lambda) g_{F}(H(\BFy,\BFxi))) \nonumber\\
&\geqslant& \mu(g_{F}(\lambda H(\BFx,\BFxi) + (1-\lambda) H(\BFy,\BFxi))) \nonumber \\
&\geqslant& \mu(g_{F}(H(\lambda \BFx + (1-\lambda) \BFy, \BFxi))), \nonumber
\end{eqnarray*}
where the first line follows from the convexity of $\mu(\cdot)$; the second line follows from the convexity of $g_{F}(\cdot)$ and the monotonicity of $\mu(\cdot)$; the third line follows from the convexity of $H(\cdot,\BFxi)$ and the monotonicity of $g_{F}(\cdot)$ and $\mu(\cdot)$.$\hfill\Box$
\endproof

Now we turn to the distribution $\mathbb{P}_1$ over $\Omega_1$. We make the following assumption in our discussion:

\begin{assumption}\label{assumption1}
$\Omega_1$ is parameterized by a finite number of parameters.
\end{assumption}

In fact, Assumption \ref{assumption1} does not cause much loss of generality. Many distribution families we are interested in are parameterized by a finite number of parameters. For example, if $\mathcal{F}$ is the family of discrete distributions whose probability mass function is:
$\mathbb{P} (\BFxi = \BFxi_i) = p_i, i=1,...,m$,
then
$$\Omega_1 = \left\{(p_1,p_2,\dots,p_m)\in \mathbb{R}^m \left| \sum_{i=1}^{m}p_i = 1; p_{i} \geqslant 0, \forall i=1,2,\dots,m  \right. \right\}$$
For family of multivariate normal distributions:
$$dF(\BFxi) = \frac{1}{\sqrt{(2\pi)^{s}|\Sigma|}}\exp(-\frac{1}{2}(\BFxi-\BFmu)^{T}\Sigma^{-1}(\BFxi-\BFmu)),$$
we have
$$\Omega_1 = \left\{(\BFmu,\Sigma)|\BFmu \in \mathbb{R}^s, \Sigma \in \mathcal{S}_{+}^{s}\right\},$$
where $\mathcal{S}_{+}^{s}$ is the cone of positive semi-definite matrices.
For family of mixture distributions $dF(\BFxi) = \sum_{i=1}^m \lambda_i p^i(\cdot):\sum_{i=1}^m \lambda_i=1,\lambda_i\geqslant 0,i=1,2,...,m$, where $p^i(\cdot),i=1,2,...,m$, are predetermined distributions,
we have
$$\Omega_1 = \left\{(\lambda_1,\lambda_2,\dots,\lambda_m)\in \mathbb{R}^m \left| \sum_{i=1}^{m}\lambda_i = 1; \lambda_{i} \geqslant 0, \forall i=1,2,\dots,m  \right. \right\}.$$

In the framework of Bayesian statistics, observations of $\BFxi$ are treated as given rather than samples randomly drawn from an underlying distribution. Meanwhile, the parameters of distribution are handled as random variables which reflect the likelihood it takes each value given the observations. We denote the distribution parameters as $\Theta = (\theta_1,\dots,\theta_m)\in \Omega_1$. To obtain the distribution of $\Theta$, one should first specify a prior distribution $f(\Theta)$ which expresses one's belief about $\Theta$ when no observation is available. Prior distribution can be either informative or uninformative. Then, when observed data $\Xi=(\BFxi_1,...,\BFxi_N)$ are collected, one can derive the posterior distribution $p(\Theta\vert\Xi)$ using Bayes' formula:
\begin{eqnarray}
p(\Theta\vert\Xi)&=&\frac{f(\Theta)\prod_{i=1}^N L(\BFxi_i \vert \Theta)}{\int_{\Omega_1}f(\Theta)\prod_{i=1}^N L(\BFxi_i \vert \Theta)d\Theta},\label{bayes_formula}
\end{eqnarray}
where $p(\Theta\vert\Xi)$ is the posterior distribution given data $\Xi$,
$\prod_{i=1}^N L(\BFxi_i \vert \Theta)$ is the likelihood function and $\int_{\Omega_1}f(\Theta)\prod_{i=1}^N L(\BFxi_i \vert \Theta)d\Theta$ is the normalizing factor.
Note that the integration in (\ref{bayes_formula}) should be replaced by summation when a discrete distribution instead of a continuous distribution is considered. In practice, if one wants to sample from the posterior distribution in (\ref{bayes_formula}), according to the Metropolis-Hastings algorithm, one only needs to be able to compute the numerator, which is usually easy to do.
In this way, the distribution $\mathbb{P}_1$ over $\Omega_1$ can be defined and sampled from easily, and it is often the case that such definitions are data-driven.

\subsection{Relation to the Models in Section \ref{sec:review}}\label{subsec:CRM_relat}
In the following, we show that all the optimization models discussed in Section \ref{sec:review} can be viewed as special cases of our proposed composite risk measure framework.

\begin{descr}
\item[\textbf{Stochastic programming.}] In stochastic programming models, we assumed that we know the distribution $F_0$ exactly, namely, $\mathcal{F} = \{F_0\}$. Therefore, models (\ref{prob_opt_expectation}), (\ref{prob_opt_VaR}) and (\ref{prob_opt_CVaR}) can be viewed as the outer risk measure taken to be a singleton:
    $$\mu(g_{F}(\cdot)) = g_{F_0}(\cdot),$$
    and the inner risk measures are chosen to be expectation, VaR and CVaR, respectively.
\item[\textbf{Distributionally robust optimization.}] In the distributionally robust optimization models, the inner measure is chosen to be the expectation measure, while the outer measure can be viewed as the worst-case risk measure:
    \begin{equation}\textup{WC}(\BFZ) = \inf \left\{\alpha \left| \mathbb{P}(\BFZ \leqslant \alpha) = 1 \right. \right\}, \label{worst_case_risk_measure}\end{equation}
    where $\BFZ$ is a random variable defined on $(\Omega_1,\Sigma_1,\mathbb{P}_1)$ with $\Omega_1$ chosen to be the distribution set $\mathcal{F}$.

    In fact, distributionally robust model also covers the singleton-coherent risk measure model:
    $$\min_{\BFx}\ \mu(\BFY(\BFx)).$$
    \cite{coherent_risk_measure} show that coherent risk measures are closely related to worst-case risk measures. The exact relationship is given in the following theorem (we use $\ll$ to denote absolute continuity between probability measures).
    \begin{theorem}\label{coherent_as_worst_case}
    Let $\mathcal{X}$ be any linear space of random variables defined on a probability space $(\Omega, \Sigma, \mathbb{P})$.
    A functional $\rho(\cdot):\mathcal{X}\rightarrow\mathbb{R}$ is a coherent risk measure if and only if there exists a family of probability measures $\mathcal{Q}$ with $\mathbb{Q}\ll\mathbb{P}$ for all $\mathbb{Q}\in \mathcal{Q}$ such that
    $$\rho(\BFX)=\sup_{\mathbb{Q}\in\mathcal{Q}}\ \mathbb{E}_{\mathbb{Q}}(\BFX),\quad \forall \BFX\in\mathcal{X},$$
    where $\mathbb{E}_{\mathbb{Q}}(\BFX)$ denotes the expectation of the random variable $\BFX$ under the measure $\mathbb{Q}$ (as opposed to the measure of $\BFX$ itself). \label{representation_of_coherent_risk}
    \end{theorem}

    For example, for $\textup{CVaR}_{\delta}(\cdot)$, the corresponding set of distribution is
    $\mathcal{Q}=\left\{\mathbb{Q}\ll\mathbb{P} \left| {d\mathbb{Q}}/{d\mathbb{P}}\leqslant (1-\delta)^{-1}\right. \right\}$.
    The above theorem shows that coherent risk measures can be represented by the worst-case expectations taken over a set of probability distributions.
    The proof of Theorem \ref{coherent_as_worst_case} in fact predates the introduction of coherent risk measure, see \cite{Huber_1981}. And this result has been used in several recent works that study distributionally robust optimizations, see, e.g., \cite{Bertsimas-uncertainty-2009} and \cite{natarajan2009constructing}.

\item[\textbf{Robust optimization.}] By using the distribution set $\mathcal{F}$ where each $F \in \mathcal{F}$ is a distribution putting all its weight on one point $\BFxi \in \BFXi$, the distributionally robust model (\ref{distributionally_robust}) reduces to a robust optimization model (\ref{prob_opt_plain_robust}) which falls into our framework.
\item[\textbf{Worst-Case CVaR and VaR optimization.}] Comparing model (\ref{worst_case_VaR}) and (\ref{worst_case_CVaR}) to the unified model (\ref{coherent_problem_1}), we see that the corresponding inner risk measures are VaR and CVaR, respectively, while the outer risk measure is the worst-case risk measure.
\end{descr}

By choosing different combinations of outer and inner risk measures, one can come up with more optimization models. However, some of those models reduce to the models above after transformation. Some examples of these models are presented as follows.

\begin{descr}
\item[\textbf{Two-fold expectation as expectation.}] Choosing outer risk measure $\mu(\cdot)$ and inner risk measure $g_{F}(\cdot)$ both as expectations, we obtain the following optimization model:
    \begin{equation}\min_{\BFx \in \mathcal{X}}\ \mathbb{E} (\mathbb{E}_{\BFxi\sim F}[H(\BFx,\BFxi)]).\label{expect_expect}\end{equation}
    Denote the cumulative distribution function of the random term $\BFxi$ by $F_{\Theta}(\BFxi)$, we have
    \begin{eqnarray*}
    \mathbb{E} (\mathbb{E}_{\BFxi\sim F}[H(\BFx,\BFxi)])
    &=& \int_{\Omega_1}\int_{\mathbb{R}^{s}}H(\BFx,\BFxi)dF_{\Theta}(\BFxi)\ p(\Theta) d\Theta\\
    &=& \int_{\mathbb{R}^{s}}H(\BFx,\BFxi)\left[\int_{\Omega_1}p(\Theta)dF_{\Theta}(\BFxi)d\Theta\right] \\
    &=& \int_{\mathbb{R}^{s}}H(\BFx,\BFxi)d\hat{F}_{\Theta}(\BFxi) \\
    &=& \hat{\mathbb{E}}[H(\BFx,\BFxi)],
    \end{eqnarray*}
    where $d\hat{F}_\Theta = \int_{\Omega_1}p(\Theta)dF_{\Theta}d\Theta$ can be viewed as the expected measure (e.g., if $\Omega_1$ contains only continuous distributions, then this is just a weighted average over all the density functions).
    Therefore, problem (\ref{expect_expect}) is equivalent to the expectation optimization problem:
    \begin{eqnarray*}
    \min_{\BFx \in \mathcal{X}}\ \hat{\mathbb{E}}[H(\BFx,\BFxi)].
    \end{eqnarray*}

\item[\textbf{Two-fold worst-case as worst-case.}] Choosing outer risk measure $\mu(\cdot)$ and inner risk measure $g_{F}(\cdot)$ both as worst-case risk measures, we obtain the following optimization model:
    \begin{equation} \min_{\BFx \in \mathcal{X}}\ \max_{F \in \mathcal{F}}\
    \max_{\BFxi \in \BFXi_F} H(\BFx,\BFxi) \label{WC-WC},\end{equation}
    where $\BFXi_F$ is the uncertainty set for $\BFxi$ when the distribution of $\BFxi$ is $F$.
    This problem can be reduced to model (\ref{prob_opt_plain_robust}):
    \begin{equation}\min_{\BFx \in \mathcal{X}}\
    \max_{\BFxi \in \BFXi} H(\BFx,\BFxi) \label{WC_WC_as_WC},\end{equation}
    where $\Xi=\left\{\BFxi:\exists F_0 \in \mathcal{F},\BFxi \in \Xi_{F_0}\right\}$. Example of this model can be found in \cite{Bertsimas-data-driven-2014}.

\end{descr}
All the cases discussed in this section are summarized in Table \ref{table_composite_risk_measure},
where the symbol $\sim$ means that the model is equivalent to another model labeled by the number,
the symbol $\times$ means that such combination of risk measures has not been considered yet.
\begin{table}\footnotesize
\centering
\caption{Different composites of risk measures}\label{table_composite_risk_measure}
\begin{tabular}{|c|cccc|}
\hline
\backslashbox{$\mu(\cdot)$}{$g_F(\cdot)$} & Expectation  &  CVaR  &  VaR  &  Worst-Case  
\\ \hline
Singleton      &(\ref{prob_opt_expectation})&(\ref{prob_opt_CVaR})&(\ref{prob_opt_VaR})&(\ref{prob_opt_plain_robust})\\
Expectation&$\sim(\ref{prob_opt_expectation})$&$\times$ &$\times$ &$\times$     \\
VaR       & $\times$ & $\times$ & $\times$ & $\times$ \\
Worst-Case &(\ref{distributionally_robust})&(\ref{worst_case_CVaR})&(\ref{worst_case_VaR})&$\sim$(\ref{prob_opt_plain_robust})
\\ \hline
\end{tabular}
\end{table}

\section{Constructing New Models}\label{sec:new_model}
In this section, we use our framework to propose and study a few new paradigms for decision making under uncertainty. In the following, we continue to use the notation $\Omega_1$ to denote the sample space of $F$ and $\mathbb{P}_1$ to denote the probability distribution over $\Omega_1$. As we have discussed earlier, $\mathbb{P}_1$ can be derived as a posterior distribution from a Bayesian approach.


\subsection{Minimizing VaR-Expectation}\label{subsec:new_var_exp}
When minimizing the expectation of random loss $H(\BFx,\BFxi)$ under distributional uncertainty, we can choose the outer risk measure as $\delta$-VaR to provide a probabilistic guarantee:
\begin{eqnarray}\min_{\BFx\in\mathcal{X}}&& \textup{VaR}_{\delta}(\mathbb{E}_{\BFxi \sim F}[H(\BFx,\BFxi)]),\label{VaR_expectation}\end{eqnarray}
where $\mathcal{X}$ is the feasible set of $\BFx$.
This model can be interpreted as finding a decision variable to minimize the threshold such that the chance that the expected loss exceeds the threshold is small,
or in other words, this model can be viewed as minimizing the upper bound for the one-sided $\delta$-confidence interval for the expected loss.
It could be applicable in the context where the expected value is a common criterion to evaluate the loss while there is uncertainty about the underlying distribution.

Note that model (\ref{VaR_expectation}) shares similar spirit as the distributionally robust model (\ref{distributionally_robust}). Both models are designed to deal with parameter uncertainty.
However, in distributionally robust models such as \cite{Delage-Ye-2010} and \cite{Bertsimas-data-driven-2014}, it is assumed that there exists a true underlying distribution $\bar{F}$ of $\BFxi$ and the distribution set $\mathcal{F}$ is chosen as a confidence region of $\bar{F}$ to hedge against uncertainty.
And the distribution set does not depend on the decision $\BFx$.
In contrast, in (\ref{VaR_expectation}), the distribution set (for the VaR) is dependent on $\BFx$. This is often desirable since for different $\BFx$, the objective function $g_{F}(H(\BFx, \BFxi))$ may have different properties, thus the set of unfavorable distributions may differ.
As a result, solving problem (\ref{VaR_expectation}) leads to a less conservative solution under the same robust level.
To illustrate, we denote the optimal solution, the optimal value and the corresponding distribution set of problem (\ref{distributionally_robust}) and problem (\ref{VaR_expectation}) by $(\BFx_{\textup{DR}}^{*},\gamma_{\textup{DR}}^{*},\mathcal{F}_{\textup{DR}})$
and
$(\BFx_{\textup{VaR}}^{*},\gamma_{\textup{VaR}}^{*},\hat{\mathcal{F}}_{\BFx}^{*})$ respectively,
then $(\BFx_{\textup{VaR}}^{*},\hat{\mathcal{F}}_{\BFx}^{*})$ is the optimal solution to the following problem:
\begin{eqnarray}
\min_{\BFx\in\mathcal{X},\mathcal{F}}&& \sup_{F\in {\cal F}}({\mathbb E}_{\BFxi\sim F}[H(\BFx,\BFxi)]) \label{dependence_illustrate}\\
\mbox{s.t.}&& \mathbb{P}_1(F \not\in\mathcal{F}) \leqslant 1-\delta,\nonumber
\end{eqnarray}
while $(\BFx_{\textup{DR}}^{*},\mathcal{F}_{\textup{DR}})$ is only a feasible solution.
Therefore, we have $\gamma_{\textup{VaR}}^{*} \leqslant \gamma_{\textup{DR}}^{*}$.

We make the following assumptions in this subsection.
\begin{assumption}\leavevmode \label{assumption_bounded}
\begin{enumerate}
\item The loss function $H(\BFx, \BFxi)$ is (piecewise) continuously differentiable and convex in $\BFx$ on ${\cal X}$.
\item Both $\mathbb{E}_{\BFxi\sim F}[H(\BFx,\BFxi)]$ and $\mathbb{E}_{\BFxi\sim F}[\nabla_{\BFx} H(\BFx,\BFxi)]$ can be evaluated efficiently.
\end{enumerate}
\end{assumption}

The first item in Assumption \ref{assumption_bounded} is necessary. Otherwise there is little hope to solve problem (\ref{VaR_expectation}) efficiently even when the risk measures over the random loss are dropped.
The second assumption ensures that standard optimization techniques can be employed to solve problem (\ref{VaR_expectation}). When the second assumption is violated, as in the case where $\mathbb{E}_{\BFxi\sim F}[H(\BFx,\BFxi)]$ does not have a closed-form expression and the dimension of $\BFx$ is high, one must turn to sample-based methods to evaluate $\mathbb{E}_{\BFxi\sim F}[H(\BFx,\BFxi)]$. Consequently, the size of the problem will be large. For those cases, we will propose an approximation method in the next subsection.

We rewrite problem (\ref{VaR_expectation}) as a chance-constrained problem:
\begin{eqnarray}
\min_{t,\BFx\in\mathcal{X}}&& t \label{VaR_chance}\\
\mbox{s.t.}&& \mathbb{P}_1(\mathbb{E}_{\BFxi\sim F}[H(\BFx,\BFxi)] \geqslant t)\leqslant 1-\delta. \nonumber
\end{eqnarray}
Thus, the optimal solution of problem (\ref{VaR_expectation}) $\BFx^{*}$ represents an optimal decision where the $\delta$ quantile of the distribution of $\mathbb{E}_{\BFxi \sim F}[H(\BFx,\BFxi)]$ is minimized (note that the random variable here is the distribution parameter $\Theta$).
A general approach to tackle problem (\ref{VaR_chance}) is the sample approximation approach (SAA) (see, e.g., \citealt{nemirovski2006scenario}, \citealt{luedtke2008sample}). Using Monte Carlo method to generate $N$ i.i.d. samples of distribution parameter $\Theta_i,i=1,\dots,N$ from distribution $\mathbb{P}_1$, problem (\ref{VaR_chance}) can be approximated by the following mixed integer nonlinear program (MINLP):
\begin{eqnarray}
\min_{t,\BFx\in\mathcal{X},z_{i}\in \{0,1\}}&& t \label{VaR_expe_MIP}\\
\mbox{s.t.}&& -M z_{i} + \mathbb{E}_{\BFxi\sim F_i}[H(\BFx,\BFxi)] \leqslant t,\quad i=1,2,...,N \nonumber\\
&& \sum_{i=1}^{N} z_{i} = \lfloor (1-\delta)N \rfloor \nonumber,
\end{eqnarray}
where $F_i$ is the distribution parameterized by $\Theta_{i}$, and $M$ is a large constant which should be larger than $\max_{\BFx\in\mathcal{X},i=1,\dots,N}\mathbb{E}_{\BFxi\sim F_i}[H(\BFx,\BFxi)] - \min_{\BFx\in\mathcal{X},i=1,\dots,N}\mathbb{E}_{\BFxi\sim F_i}[H(\BFx,\BFxi)]$. When $\mathcal{X}$ is bounded, for example, there exists a finite-valued $M$.

The rationale of this approximation approach is quite clear. By generating $N$ instances of $\Theta$, a uniform distribution over $\Theta_i,i=1,\dots,N$ instead of $\mathbb{P}_1$ is considered when minimizing the $\delta$-quantile of the expectation. When $N$ is sufficiently large, problem (\ref{VaR_expe_MIP}) can approximate problem (\ref{VaR_chance}) well by the law of large numbers. Specifically, we have the following accuracy and feasibility bound which provides a guideline for choosing suitable $N$ given desired accuracy levels. The next theorem largely follows the results in \cite{luedtke2008sample}.

\begin{theorem} \label{theorem_luedtke}
Suppose that $\mathcal{X}$ is bounded with diameter $D$ and $\mathbb{E}_{\BFxi\sim F}[H(\BFx,\BFxi)]$ is globally Lipschitz continuous with Lipschitz constant $L$ for all $F\in\mathcal{F}$. Assume that both problems (\ref{VaR_chance}) and (\ref{VaR_expe_MIP}) are feasible and their optimal values are finite. Denote (\ref{VaR_chance}) and (\ref{VaR_expe_MIP}) by $P_{\delta}$ and $P_{N, \delta}$, and their optimal solutions (only the $\BFx$ part) and optimal values by $(\BFx_{\delta}^{*},t_{\delta}^{*})$ and $(\BFx^*_{N, \delta},t^*_{N,\delta})$, respectively.
Let $0 < \tau \leqslant 1-\delta$, $\epsilon\in(0,1)$, $\gamma > 0$,
\begin{eqnarray}
N_0 = \frac{2}{\tau^2}\left[\log\left(\frac{1}{\epsilon}\right)+n\log \left\lceil \frac{2LD}{\gamma} \right\rceil+\log \left\lceil \frac{2}{\tau} \right\rceil\right].\end{eqnarray}
Then, when $N \geqslant N_0$, with probability at least $1-2\epsilon$, we have:
\begin{equation}t^*_{\delta-\tau} -\gamma \leqslant t^*_{N, \delta} \leqslant t^*_{\delta + \tau}.
\label{bound_of_two_sides}\end{equation}
\end{theorem}

\proof{Proof.}
When $N \geqslant N_0 \geqslant \frac{1}{2 \tau^2}\log \left(\frac{1}{\epsilon}\right)$, the upper bound in (\ref{bound_of_two_sides}) is given by Theorem 3 in \cite{luedtke2008sample}.
Also, it follows from Theorem 10 in \cite{luedtke2008sample} that a feasible $(\BFx,t)$ for problem:
\begin{eqnarray}
\min_{t,\BFx,z_{i}\in \{0,1\}}&& t \label{VaR_chance_discrete}\\
\mbox{s.t.}&& -M z_{i} + \mathbb{E}_{\BFxi\sim F_i}[H(\BFx,\BFxi)]+\gamma \leqslant t,\quad i=1,2,...,N \nonumber\\
&& \sum_{i=1}^{N} z_{i} = \lfloor (1-\delta)N \rfloor \nonumber
\end{eqnarray}
is feasible for problem $P_{\delta-\tau}$ with at least $1-\epsilon$ probability. Notice that the optimal $t$ for problem (\ref{VaR_chance_discrete}) equals $t^*_{N, \delta}+\gamma$. Thus, when the feasible solution $(\BFx,t)$ for problem (\ref{VaR_chance_discrete}) is feasible for problem $P_{\delta-\tau}$, we have
$t^*_{N, \delta}+\gamma \geqslant t^*_{\delta-\tau}$. This completes the proof.
$\hfill\Box$
\endproof

Note that in the above conditiopns, the logarithms are taken over $\epsilon$ and $\gamma$, therefore we can choose very small values of them without increasing the number of necessary samples very much.
Meanwhile, experiments show that the bound in Theorem \ref{theorem_luedtke} is very conservative (\citealt{luedtke2008sample}) and one can choose a much smaller $N$ in practice.

Numerically, software packages such as CPLEX (mixed integer linear problems) and MOSEK (certain MINLP problems) can be used to solve moderate sized problems. When Assumption $\ref{assumption_bounded}$ holds, problem (\ref{VaR_expe_MIP}) may be solved by those solvers.

\begin{remark}
It is worth pointing out that model (\ref{VaR_expe_MIP}) can also be used to solve a family of VaR-VaR problems and VaR-CVaR problems. For portfolio optimization problems where $H(\BFx,\BFxi)=\BFxi^{T}\BFx$ and the distribution of $\BFxi$ is a normal distribution with mean $\BFmu$ and covariance $\Gamma$,
both $\textup{VaR}_{\delta}(H(\BFx,\BFxi))$ and $\textup{CVaR}_{\delta}(H(\BFx,\BFxi))$ have closed-form expressions (\citealt{sarykalin2008value}):
\begin{eqnarray}
\textup{VaR}_{\delta}(\BFxi^{T}\BFx) &=& \Phi^{-1}(\delta)\sqrt{\BFx^{T}\Gamma\BFx}+\BFmu^{T}\BFx \label{VaR_normal},\\
\textup{CVaR}_{\delta}(\BFxi^{T}\BFx) &=& \frac{1}{\sqrt{2\pi}(1-\delta)}\exp(-(\Phi^{-1}(\delta))^2/2)\sqrt{\BFx^{T}\Gamma\BFx}
+\BFmu^{T}\BFx\label{CVaR_normal},
\end{eqnarray}
where $\Phi^{-1}(\cdot)$ is the inverse of the cumulative distribution function of a standard normal distribution.
In this case,
we can minimize VaR$_{\delta}$-VaR$_{\epsilon}$ and VaR$_{\delta}$-CVaR$_{\epsilon}$ by solving the following problems:
\begin{eqnarray}
\min_{\BFx\in\mathcal{X}}&& \textup{VaR}_{\delta}\left(\Phi^{-1}(\epsilon)\sqrt{\BFx^{T}\Gamma\BFx}+\BFmu^{T}\BFx\right),\label{VaR_VaR_normal}\\
\min_{\BFx\in\mathcal{X}}&& \textup{VaR}_{\delta}\left(\frac{1}{\sqrt{2\pi}(1-\epsilon)}\exp(-(\Phi^{-1}(\epsilon))^2/2)
\sqrt{\BFx^{T}\Gamma\BFx}+\BFmu^{T}\BFx\right),\label{VaR_CVaR_normal}
\end{eqnarray}
both of which can be solved approximately by mixed integer second order cone program (MISOCP) using SAA approach.
\end{remark}
\begin{remark}
It is also worth mentioning that when the objective function
$H(\BFx,\BFxi)$ is linear in $\BFxi$ (such as in the portfolio
selection problem), the VaR-Expectation model is similar to a VaR
model in formulation. More precisely, suppose $H(\BFx, \BFxi) =
\BFxi^T f(\BFx)$, then $\textup{VaR}_{\delta}({\mathbb E}_{\BFxi}
[H(\BFx,\BFxi)]) = \textup{VaR}_{\delta}(({\mathbb E}\BFxi)^T f(\BFx))$. And
the latter can be viewed as a VaR model
with ${\mathbb E}\BFxi$ as the random variable. Recently, such a formulation
has been studied in \cite{cui2013nonlinear} and \cite{wen2013asset} in the context
of the portfolio selection problem.
In particular, \cite{wen2013asset} show that an alternating direction augmented Lagrangian method (ADM) can be applied to solve the problem very efficiently in that case.
However, the VaR-Expectation model is
different from the VaR model if the objective is not linear in
$\BFxi$ since the expectation can no longer be taken into $\BFxi$.
Moreover, even in the linear case, the VaR-Expectation model has a
very different interpretation from the VaR model: In the VaR model,
the uncertainty is directly in $\BFxi$, and one has to assume a
certain distribution of $\BFxi$ (which is usually estimated from an
empirical distribution); while in the VaR-Expectation model, the uncertainty is
in the distribution of $\BFxi$. Such different interpretations would
usually lead to different uncertainty sets even with the same set of
observations.
\end{remark}

\subsection{Minimizing CVaR-CVaR and CVaR-Expectation}\label{subsec:new_cvar}
Using the same idea as above, we formulate a robust model for CVaR optimization by choosing $\mu(\cdot)$ as $\delta$-VaR and $g_{F}(\cdot)$ as $\epsilon$-CVaR, namely,
\begin{eqnarray*}
&&\min_{\BFx \in \mathcal{X},t\in\mathbb{R}}\ t \label{VaR_CVaR}\\
&&\mbox{s.t.}\quad\mathbb{P}_1 (\textup{CVaR}_{\epsilon}(H(\BFx,\BFxi_{F}))\geqslant t)\leqslant 1-\delta, \nonumber
\end{eqnarray*}
where $\mathbb{P}_1(\cdot)$ is the probability measure of $F$, and $\BFxi_{F}$ means that the distribution of $\BFxi$ is $F$.
This model can be viewed as minimizing the upper bound for the one-sided $\delta$-confidence interval for the CVaR.

For most cases, it is impossible to derive a closed-form expression for CVaR. In practice, sample based methods like sample average approximation (SAA) are widely used to compute CVaR. To ensure the accuracy of evaluation, the number of samples is typically large (the exact number of necessary samples will be discussed later). Therefore, if we directly replace expectation by CVaR in model (\ref{VaR_expe_MIP}), the resulting problem will be a large-sized mixed integer program, which is difficult to be solved. The same situation occurs in a VaR-Expectation model where the expectation cannot be evaluated directly and must be approximated by sample average.

To deal with the issue mentioned above, we relax the outer VaR to CVaR, leading to a CVaR-CVaR (nested CVaR) model:
\begin{eqnarray}
\min_{\BFx \in \mathbb{R}^n, \alpha \in \mathbb{R}} \alpha +\frac{1}{1-\delta}
\mathbb{E}[(\textup{CVaR}_{\epsilon}(H(\BFx,\BFxi_{F})) - \alpha)^{+}].\label{CVaR_CVaR}
\end{eqnarray}
Similarly, we can consider the following CVaR-Expectation model instead of the VaR-Expectation model:
\begin{eqnarray}\min_{\BFx\in\mathcal{X}}&& \textup{CVaR}_{\delta}(\mathbb{E}_{\BFxi \sim F}[H(\BFx,\BFxi)]).\label{CVaR_expectation}\end{eqnarray}

The reason for making the above relaxations is two fold. First,
\cite{coherent_general} shows that $\textup{CVaR}_{\delta}(Z)$ is the smallest upper bound of $\textup{VaR}_{\delta}(Z)$ among all coherent risk measures that depend
only on the distribution of $Z$. By relaxing the original problem to convex programs, a wide range of convex optimization algorithms can be employed to solve them. Second, optimization problems similar to model (\ref{CVaR_CVaR}) have been investigated in the literature in the context of risk averse multistage stochastic programming (see, e.g., \citealt{guigues2012sampling}, and \citealt{philpott2012dynamic}), thus similar algorithms and techniques can be employed.
In addition, since model (\ref{CVaR_CVaR}) and (\ref{CVaR_expectation}) take into account the extent of losses in the most adverse scenarios, they are also meaningful in their own right. In the following, we present a general regime to solve problem (\ref{CVaR_CVaR}) using SAA approach.
Since the discussion of model (\ref{CVaR_expectation}) resembles that of model (\ref{CVaR_CVaR}) and is simpler, we will only discuss model (\ref{CVaR_CVaR}).

SAA is a popular method in stochastic programming and has been used to deal with CVaR optimization problems, see \cite{shapiro2006complexity} and \cite{wang2008sample}. The concept of SAA is to approximate an expectation by the average of many samples generated by Monte Carlo method or other schemes. In problem (\ref{CVaR_CVaR}), we have
\begin{eqnarray*}
\mathbb{E}[(\textup{CVaR}_{\epsilon}(H(\BFx,\BFxi)) - \alpha)^{+}]
\approx \frac{1}{N}\sum_{i=1}^{N}(\textup{CVaR}_{\epsilon}(H(\BFx,\BFxi_{F_i}))-\alpha)^{+},
\end{eqnarray*}
where $N$ i.i.d. samples of $\Theta$: $\Theta_1,\dots,\Theta_N$ are generated from the distribution of $\Theta$ and the parameter of the distribution of $\BFxi_{F_i}$ is $\Theta_i$.
Then, the original optimization problem can be approximated by the following problem:
\begin{eqnarray}
\min_{\alpha \in \mathbb{R},\BFx \in \mathbb{R}^n, \BFu \in \mathbb{R}^N}
&&\alpha + \frac{1}{(1-\delta)N} \BFe^T \BFu \label{CVaR_CVaR_relax1}\\
\mbox{s.t.}\quad  &&u_i \geqslant \textup{CVaR}_{\epsilon}(H(\BFx,\BFxi_{F_i}))-\alpha,\quad i=1,2,...,N \nonumber\\
&& u_i \geqslant 0,\quad i=1,2,...,N, \nonumber
\end{eqnarray}
where $\BFe$ denotes a vector of all ones. Notice that problem (\ref{CVaR_CVaR_relax1}) is a CVaR constrained problem, and the CVaR in the constraint can also be approximated using SAA, resulting in the following problem:
\begin{eqnarray}
\min_{\alpha,\BFx, \BFu, \BFv}
&&\alpha + \frac{1}{(1-\delta)N} \BFe^T \BFu \label{CVaR_CVaR_relax2}\\
\mbox{s.t.}\quad  &&u_i \geqslant v_i + \frac{1}{(1-\epsilon)M} \sum_{j=1}^M [H(\BFx ,\BFxi_{F_i}^j) - v_i]^+ - \alpha,\quad i=1,2,...,N \nonumber\\
&& u_i \geqslant 0,\quad i=1,2,...,N, \nonumber
\end{eqnarray}
where $\BFv = (v_{1},\dots,v_{N}),i=1,\dots,N$ are auxiliary variables and $(\BFxi_{F_i}^{1},\dots,\BFxi_{F_i}^{M}),i=1,\dots,N$ are i.i.d. samples generated from $F_i$.
When $H(\BFx,\BFxi)$ is convex in $\BFx$, problem (\ref{CVaR_CVaR_relax2}) is a convex program with a linear objective function and can be solved efficiently for large sized problems. Now we turn to the number of samples, $N$ and $M$. For $M$, \cite{wang2008sample} shows that under certain regularity conditions, there is at least $1-\epsilon$ probability that the SAA of CVaR lies within the $\gamma$-neighborhood of CVaR when
\begin{eqnarray}
M \geqslant \frac{C_1(H,F)}{\gamma^{2}}\left[C_2(H,F)n + C_3(H,F) \log\left(\frac{1}{\epsilon}\right)\right] \label{bound_of_M},
\end{eqnarray}
where $C_1(H,F)$, $C_2(H,F)$ and $C_3(H,F)$ are constants for a given objective function $H(\BFx,\BFxi)$ and distribution $F$, $n$ is the dimension of $\BFx$. Since the constants are typically difficult to calculate, the bound in (\ref{bound_of_M}) serves as a benchmark to estimate the order of $M$. For $N$, \cite{shapiro2006complexity} proves that under mild conditions, the sample size
\begin{eqnarray}
N \geqslant \frac{D_1(H,\mathbb{P}_1)}{\gamma^{2}}\left[n  \log \left(\frac{D_2(H,\mathbb{P}_1)}{\gamma}\right) + \log \left(\frac{D_3(H,\mathbb{P}_1)}{\epsilon}\right)\right] \label{bound_of_N}
\end{eqnarray}
ensures that the solution of problem (\ref{CVaR_CVaR_relax1}) lies within the $\gamma$-neighborhood of the solution of problem (\ref{CVaR_CVaR}) with probability $1-\epsilon$. Here $D_1(H,\mathbb{P}_1)$, $D_2(H,\mathbb{P}_1)$ and $D_3(H,\mathbb{P}_1)$ are constants for a given objective function $H(\BFx,\BFxi)$ and distribution $\mathbb{P}_1$ (of $F$), $n$ is the dimension of $\BFx$.

Using the estimates in (\ref{bound_of_M}) and (\ref{bound_of_N}), we can see that the size of problem (\ref{CVaR_CVaR_relax2}) is quite large. In practice, there are many ways to accelerate the computation. For instance, \cite{alexander2006minimizing} and \cite{xu2009smooth} show that smoothing the CVaR can enhance the computation efficiency up to several times.

\begin{remark}
When the expression inside CVaR can be evaluated and differentiated easily, we only need to sample from distribution $\mathbb{P}_1$ and the resulting problem is much easier to solve. For instance, consider the following models:
\begin{eqnarray}
\min_{\BFx\in\mathcal{X}}&& \textup{CVaR}_{\delta}(\textup{VaR}_{\epsilon}(\BFxi^{T}\BFx)),\label{CVaR_VaR}\\
\min_{\BFx\in\mathcal{X}}&& \textup{CVaR}_{\delta}(\textup{CVaR}_{\epsilon}(\BFxi^{T}\BFx)),\label{CVaR_CVaR_normal}
\end{eqnarray}
where the distribution of $\BFxi$ is normal. Since $\textup{VaR}_{\epsilon}(\BFxi^{T}\BFx)$ and $\textup{CVaR}_{\epsilon}(\BFxi^{T}\BFx)$ have closed-form expressions (given in (\ref{VaR_normal}), (\ref{CVaR_normal})), problem (\ref{CVaR_VaR}) and (\ref{CVaR_CVaR_normal}) can be solved by second order cone programs (SOCP). For CVaR-Expectation problems, when the distribution of $\BFxi$ is discrete, problem (\ref{CVaR_expectation}) reduces to a mean-CVaR problem and can be solved efficiently (\citealt{iyengar2013fast} and \citealt{wen2013asset}).
\end{remark}

Having filled in several empty entries in Table \ref{table_composite_risk_measure},
we have Table \ref{table_composite_risk_measure_full}. In Table \ref{table_composite_risk_measure_full}, the symbol $\sim$ means that the
model is equivalent to another model labeled by the number and
the symbol $\times$ means that such combination of risk measures has not been considered yet.
(In Table \ref{table_composite_risk_measure_full}, Var-Worst-Case and CVaR-Worst-Case model can be viewed as VaR-CVaR and CVaR-CVaR models, respectively, with the $\delta$ in the inner CVaR being 1.)

\begin{table}\footnotesize
\centering
\caption{Different composites of risk measures: a full table}\label{table_composite_risk_measure_full}
\begin{tabular}{|c|cccc|}
\hline\up
\backslashbox{$\mu(\cdot)$}{$g_F(\cdot)$} & Expectation      &  CVaR     &VaR&  Worst-Case   \down \\ \hline\up
Singleton       &(\ref{prob_opt_expectation})&(\ref{prob_opt_CVaR})&(\ref{prob_opt_VaR})&(\ref{prob_opt_plain_robust})\\
Expectation&$\sim(\ref{prob_opt_expectation})$&$\times$ &$\times$ &$\times$     \\
VaR       &(\ref{VaR_expectation})&(\ref{VaR_CVaR_normal})&(\ref{VaR_VaR_normal})  &$\sim$(\ref{VaR_CVaR_normal})\\
Worst-Case &(\ref{distributionally_robust})&(\ref{worst_case_CVaR})&(\ref{worst_case_VaR})&$\sim$(\ref{prob_opt_plain_robust})\\
CVaR &(\ref{CVaR_expectation})&(\ref{CVaR_CVaR})&(\ref{CVaR_VaR})&$\sim$(\ref{CVaR_CVaR})\down\\
                                 \hline
\end{tabular}
\end{table}

\section{Numerical Experiments}\label{sec:numer}
In this section, we conduct numerical experiments to demonstrate the
tractability and effectiveness of the models proposed in Section \ref{sec:new_model}.
In particular, we consider the VaR-Expectation model, the
CVaR-Expectation model and the CVaR-CVaR model and conduct
numerical experiments using the portfolio selection problem. In the
portfolio selection problem, a decision maker has to choose a
portfolio to invest using available stocks based on historical data.
In particular, our data set contains the daily returns of 359
different stocks that are in the S\&P 500 index and do not have
missing data from 2010 to 2011. In the following, we first show that
when using a VaR-Expectation model, the corresponding distribution
sets (for the VaR) depend on the decision, which is a distinguishing
feature of the model. We show that this feature makes the
VaR-Expectation model less conservative than the
distributionally robust model in \citet{Delage-Ye-2010}. Then we
demonstrate that all the three models can be computed efficiently even when
the sample size is large enough to ensure the precision of the solution.
Finally, we compare the resulting returns of the three models with
existing models.

\subsection{Comparing the VaR-Expectation model and the distributionally robust model}\label{subsec:numer_compare}
\begin{figure}
\FIGURE
{\includegraphics*[width=\textwidth]{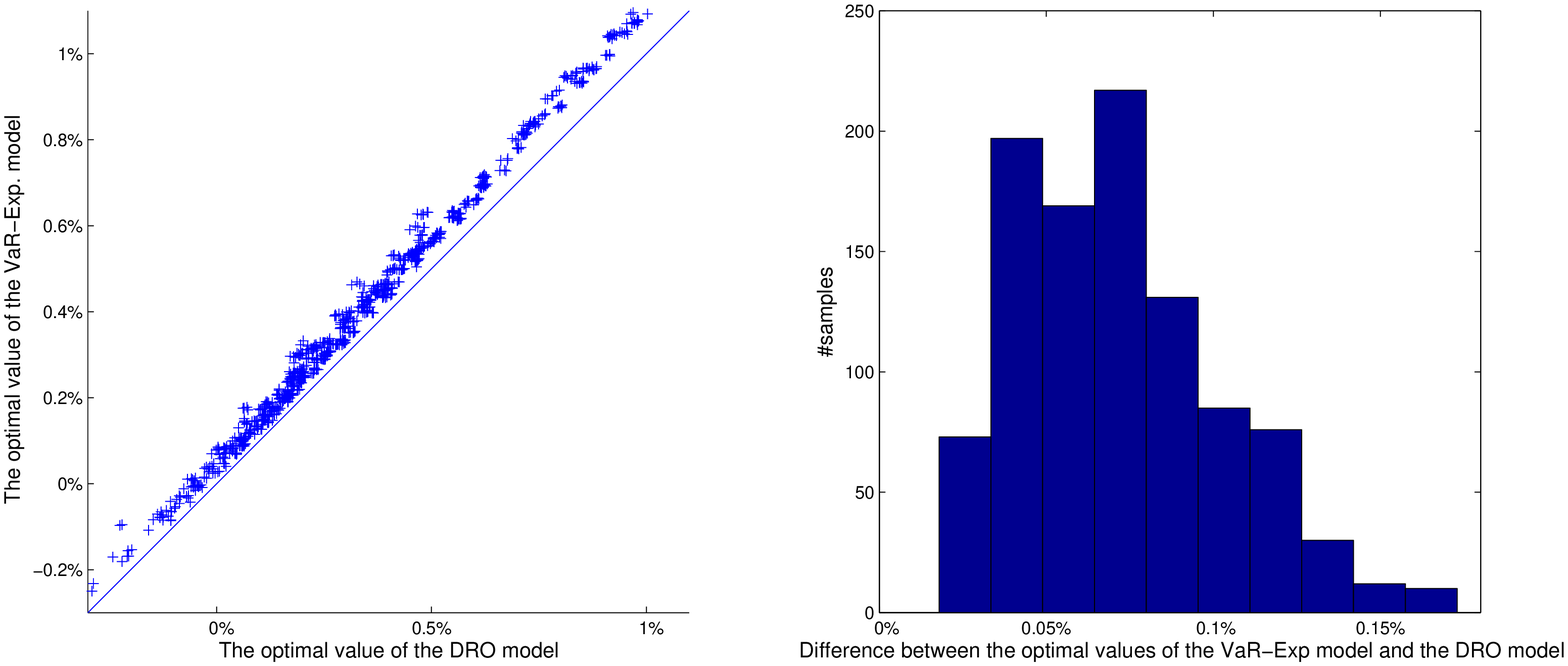}}
{Comparison between the VaR-Expectation model and the DRO model.\label{fig:qqplot}}
{} 
\end{figure}

\begin{figure}
\FIGURE
{\includegraphics*[width=\textwidth]{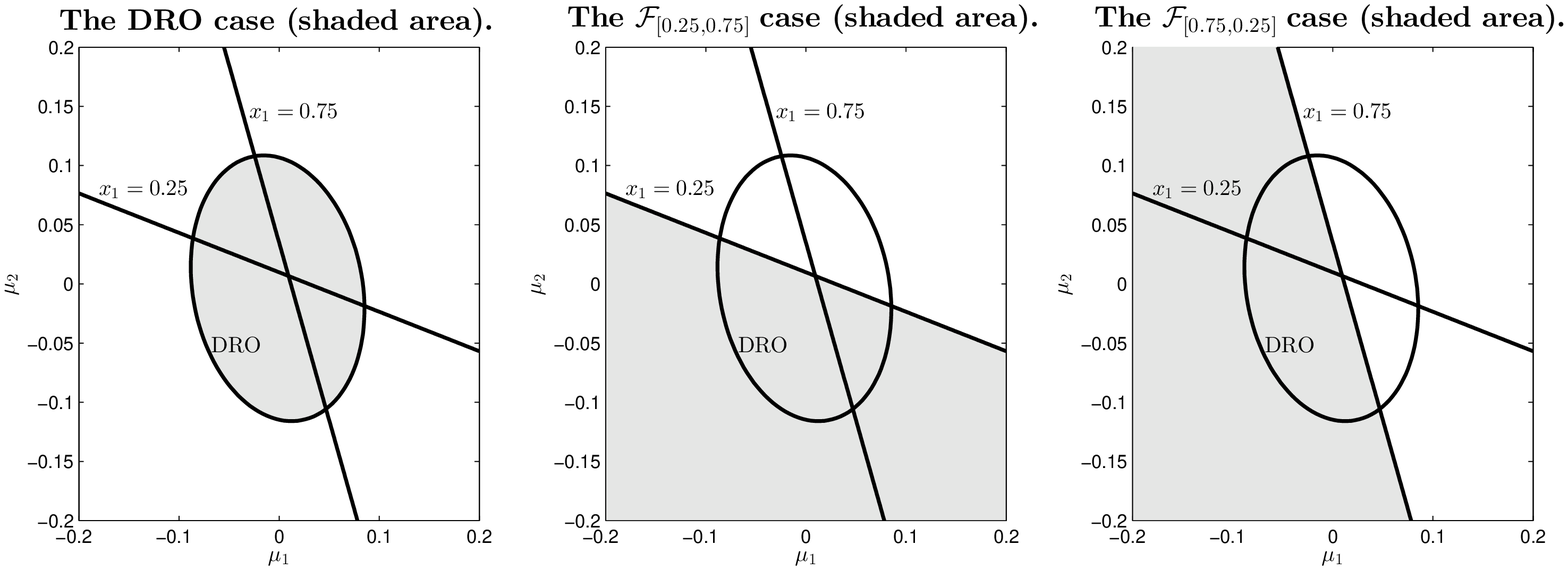}}
{Comparison of distribution sets.\label{fig:depend}}
{} 
\end{figure}

We have discussed in Section \ref{sec:new_model} that one important feature of our
proposed approach, the VaR-Expectation model, is that the
distribution set (for the VaR) depends on the decision $\BFx$
and therefore the obtained solution from this model will be less
conservative than that of a traditional distributionally
robust model in which the distribution set is independent of $\BFx$.
Here we illustrate this feature using numerical examples. Suppose a decision
maker builds a portfolio of $n$ stocks using the VaR-Expectation model.
The joint distribution of the stocks is parameterized by an
$n$-dimensional normal distribution with mean $\BFmu$ and covariance
matrix $\Sigma$. At each decision time, we use the returns in the
past $t$ days to derive the posterior distribution for $\BFmu$ and
$\Sigma$. Let $\BFmu_0$ and $\Sigma_0$ denote the mean and
covariance of the empirical distribution of the stocks in the past
$t$ days. Using multivariate Jeffereys prior density as the prior
distribution (which is the uninformative prior for normal
distributions), we have that the posterior distribution $\mathbb{P}_1$ of $\BFmu$ and
$\Sigma$ is given by (\citealt{gelman2013bayesian}):
\begin{equation}f(\BFmu,\Sigma) = \mathcal{N}(\BFmu \vert \BFmu_0,t^{-1}\Sigma) \mathcal{W}^{-1}(\Sigma \vert t\Sigma_0,t-1),\label{normal_inverse_wishart}\end{equation}
where $\mathcal{N}(\BFmu \vert \BFmu_0,t^{-1}\Sigma)$ is the
probability density function of a multivariate normal distribution
with mean $\BFmu_0$ and covariance matrix $t^{-1}\Sigma$ and
$\mathcal{W}^{-1}(\Sigma \vert t\Sigma_0,t-1)$ is the probability
density function of an inverse-Wishart distribution with scale
matrix $t\Sigma_0$ and degree of freedom $t-1$. In the following,
we set $t=30$.
For each fixed decision (allocation) $\BFx$, the worst-case
distribution
$F_0(\BFx)=\mathcal{N}(\BFmu_0(\BFx),\Sigma_0(\BFx))$
of loss $\BFxi$ satisfies the following condition:
\begin{equation}
\mathbb{P}_1(\mathbb{E}_{F}[\BFxi^{T}\BFx] \geqslant
\mathbb{E}_{F_0(\BFx)}[\BFxi^{T}\BFx])=1-\delta,
\nonumber\end{equation} or equivalently,
\begin{equation}
\mathbb{P}_1(\BFmu_F^{T}\BFx \geqslant
\BFmu_0(\BFx)^{T}\BFx)=1-\delta. \nonumber\end{equation}
Thus, the distribution set of the VaR-Expectation model is:
\begin{equation}
\mathcal{F}_{\BFx} = \left\{\mathcal{N}(\BFmu,\Sigma)|\BFmu^{T}\BFx
\leqslant \BFmu_0(\BFx)^{T}\BFx\right\}, \nonumber\end{equation} which
depends on $\BFx$. The optimization problem under the
VaR-Expectation model is (we assume that the total investment amount
must equal to $1$):
\begin{eqnarray}\label{numerical_varexpectation}
\max_{\BFx \geqslant 0, \mathbf{e}^T\BFx = 1}&\BFmu_0(\BFx)^T \BFx.
\end{eqnarray}
To compare to the distributionally robust model in
\cite{Delage-Ye-2010}, we note that the objective function in
\cite{Delage-Ye-2010} is
\begin{eqnarray}\label{numerical_dro}
\max_{\BFx \geqslant 0, \mathbf{e}^T\BFx = 1} \min_{F\in {\mathcal F}}&\BFmu_F^T \BFx
\end{eqnarray}
where
\begin{equation}
    \mathcal{F}= \left\{F(\BFxi)\in \mathcal{M}\left|
    \begin{aligned}
    &(\mathbb{E}_{\BFxi\sim F}(\BFxi)-\BFmu_0)^T \Sigma_0^{-1} (\mathbb{E}_{\BFxi\sim F}(\BFxi)-\BFmu_0)\leqslant \gamma_1 \\
    &\mathbb{E}_{\BFxi\sim F} [(\BFxi-\BFmu_0)(\BFxi-\BFmu_0)^T] \preceq \gamma_2 \Sigma_0
    \end{aligned}
    \right.
    \right\}, \nonumber
    \end{equation}
where $\gamma_1$ and $\gamma_2$ are chosen according to the
discussions in \cite{Delage-Ye-2010} to ensure that with probability
$\delta$, ${\mathcal F}$ contains the true distribution. In the
following, we fix $\delta = 0.95$. We perform the following procedures:
we randomly pick two stocks and use their empirical returns
from 3/4/2010 to 4/14/2010 to fit a normal distribution.
Then we draw $10^6$ data from the distribution to form
the data set.
(In \citealt{Delage-Ye-2010}, it usually needs at least $10^5$ data points to get a valid $\gamma_1$ and $\gamma_2$, therefore we have to take this bootstrapping method.)
Then we solve (\ref{numerical_varexpectation}) and
(\ref{numerical_dro}) respectively using these data. By the
discussions in Section \ref{subsec:new_var_exp}, the optimal value of
(\ref{numerical_varexpectation}) should be larger than that of
(\ref{numerical_dro}). In our numerical experiment, we repeat the
above procedures $1000$ times and plot the result in Figure
\ref{fig:qqplot}.
The left figure in Figure \ref{fig:qqplot} shows a scatter plot of the optimal values obtained by the
distributionally robust model in \cite{Delage-Ye-2010} ($x$-axis) and the
VaR-Expectation model ($y$-axis) in the 1000 experiments, and the right figure shows the distribution of the difference of the optimal values in the same set of experiments.
From Figure \ref{fig:qqplot}, we can see that the VaR-Expectation model results in higher value in all the cases, meaning that it is less conservative. Indeed, the average optimal value of
(\ref{numerical_varexpectation}) in the $1000$ experiments is 0.070\% larger than that of (\ref{numerical_dro}).
Moreover, for one particular experiment, we plot the projection of
$\mathcal{F}_{[0.75,0.25]}$, $\mathcal{F}_{[0.25,0.75]}$ and the
distribution set of the DRO model on the plane spanned by $\BFmu$ in
Figure \ref{fig:depend}. From Figure \ref{fig:depend}, we can see
that in the DRO approach, the distribution set is independent of
the choice of $\BFx$. However, in the VaR-Expectation model, the
distribution set changes with the choice of $\BFx$. As we have
mentioned above, it is such feature of the VaR-Expectation model
that makes the solution less conservative yet with the same level of
probabilistic guarantee.

\subsection{Solving the composite risk measure model}\label{subsec:numer_solve}
Under the setting in Section \ref{subsec:numer_compare}, the VaR-Expectation model can be
solved by the alternating direction augmented Lagrangian method
(ADM, see \citealt{wen2013asset}), the CVaR-Expectation model can be
solved by an LP, and the CVaR-CVaR model can be
solved by an SOCP. We evaluate the
performance of these three models with different stock number $n$
and sample size $N$. Our experiments are performed on a laptop with
8.00 GB of RAM and 2.20 GHz processor, using MOSEK with the Matlab
interface.

\begin{figure}
\FIGURE
{\includegraphics*[width=\textwidth]{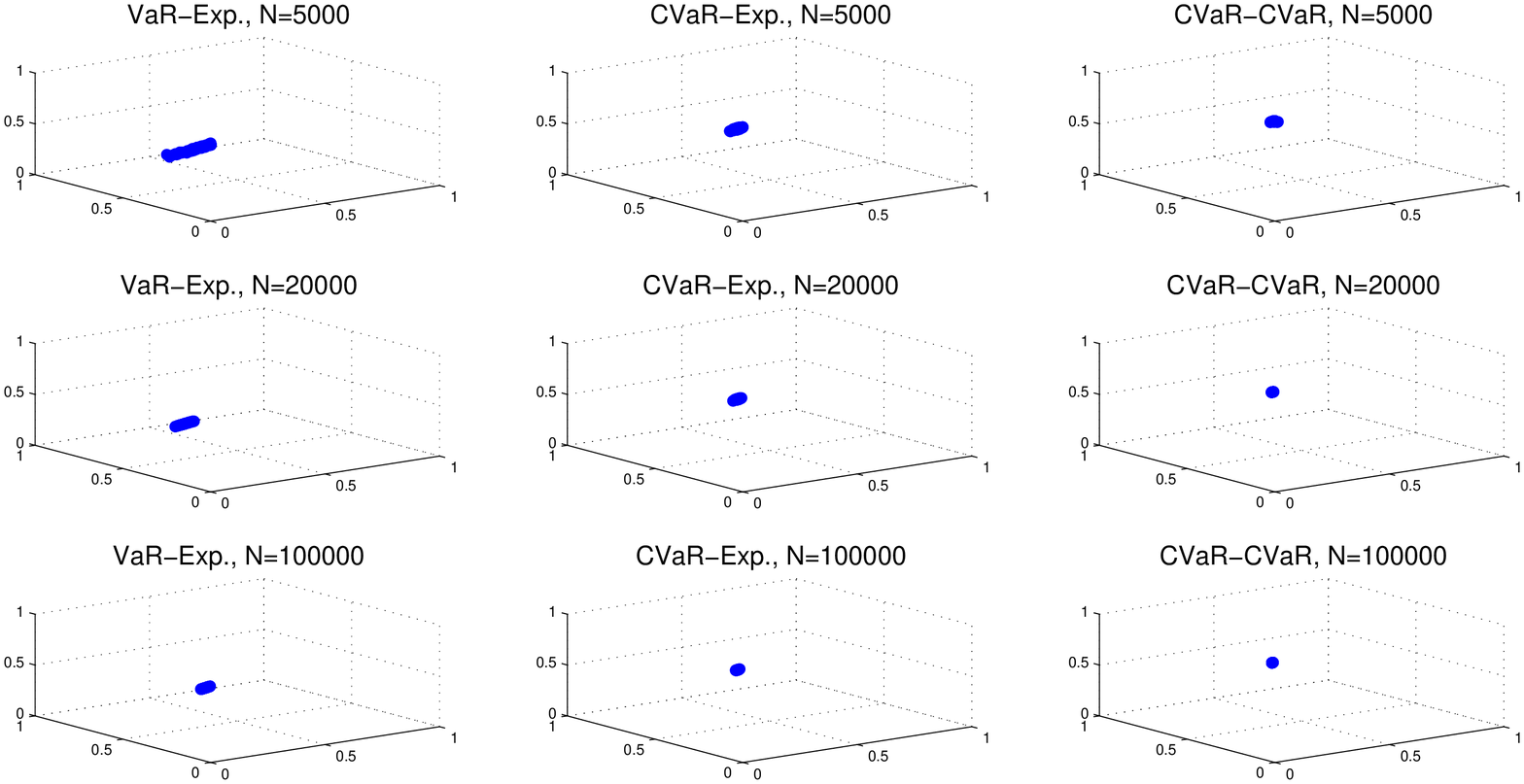}}
{Optimal allocations of the approximated problems.\label{fig:allocation}}
{} 
\end{figure}

We first solve the VaR-Expectation, CVaR-Expectation and CVaR-CVaR models for different sample size $N$ when $t=30$ and $n=4$, while both outer and inner risk levels are chosen as $0.95$.
The stocks are randomly chosen from all 359 stocks and the period we consider is from 3/4/2010 to 4/14/2010.
We perform the experiment on the same set of stocks 100 times, and the results are shown in Table \ref{experiment_table1} and Figure \ref{fig:allocation}. In Table \ref{experiment_table1}, the first column, denoted by ``ave'', for each method is the average of the optimal values of the corresponding models, the second column, denoted by ``std'', is the standard deviation of the optimal values in all experiments, while the third column is the average computation time (in seconds). Notice that when $N=100000$, all these three models can still be computed efficiently. In the mean time, the solutions of SAA problems can approximate the solution of the original models very well.
Figure \ref{fig:allocation} displays the scatter plots of the weights of the first 3 stocks in the optimal portfolio of all experiments (the weight of the last stock can be computed by one minus the total weights of the first three stocks). The result shows that the larger $N$ is, the more concentrated the solution is. It indicates that the optimal solution also converges as the sample size becomes large.

\begin{table}\footnotesize \caption{Computation time and precision of different CRM models when $n=4$.}
\label{experiment_table1}


\centering
\begin{tabular}{ccccccccccccc}
\hline
\multicolumn{1}{c}{\#samples}&&
\multicolumn{3}{c}{VaR-Expectation}&&
\multicolumn{3}{c}{CVaR-Expectation}&&
\multicolumn{3}{c}{CVaR-CVaR}\\
\cline{3-5}\cline{7-9}\cline{11-13}
$N$ && ave & std & time &&
ave & std & time &&
ave & std & time\\
\hline
5000 &&   0.0005 & $ 1.91 \times 10^{-4}$ &     5.12 &&  -0.0007 & $ 1.44 \times 10^{-4}$ &     0.13 &&  -0.0370 & $ 2.80 \times 10^{-4}$ &     0.76 \\ \hline
10000 &&   0.0004 & $ 1.56 \times 10^{-4}$ &     5.17 &&  -0.0007 & $ 7.94 \times 10^{-5}$ &     0.16 &&  -0.0371 & $ 2.01 \times 10^{-4}$ &     1.78 \\ \hline
20000 &&   0.0004 & $ 1.09 \times 10^{-4}$ &     4.99 &&  -0.0007 & $ 6.68 \times 10^{-5}$ &     0.25 &&  -0.0371 & $ 1.44 \times 10^{-4}$ &     3.26 \\ \hline
50000 &&   0.0005 & $ 6.43 \times 10^{-5}$ &     7.58 &&  -0.0007 & $ 4.20 \times 10^{-5}$ &     0.70 &&  -0.0372 & $ 7.94 \times 10^{-5}$ &    10.02 \\ \hline
100000 &&   0.0004 & $ 4.63 \times 10^{-5}$ &     9.26 &&  -0.0007 & $ 3.18 \times 10^{-5}$ &     1.65 &&  -0.0373 & $ 7.17 \times 10^{-5}$ &    17.64 \\ \hline
\end{tabular}
  \end{table}

Now we fix $N=5000$ and choose different $n$, while the period considered and all other parameters are the same as in the previous experiment. The results are presented in Table \ref{experiment_table2}, where $t_{ave}$ is the average computation time in all experiments, $t_{min}$ is the minimum computation time and $t_{max}$ is the maximum computation time. The models considered are (from left to right in Table \ref{experiment_table2}): VaR-Expectation model, CVaR-Expectation model and CVaR-CVaR model. We can see that our new models can still be computed in a reasonably amount of time even when $n$ is large.

\begin{table}\footnotesize \caption{Computation time of different CRM models when $N=5000$.}
\label{experiment_table2}


\centering
\begin{tabular}{ccccccccccccc}
\hline
\multicolumn{1}{c}{\#stocks}&&
\multicolumn{3}{c}{VaR-Expectation}&&
\multicolumn{3}{c}{CVaR-Expectation}&&
\multicolumn{3}{c}{CVaR-CVaR}\\
\cline{3-5}\cline{7-9}\cline{11-13}
$n$ && $t_{ave}$ & $t_{min}$ & $t_{max}$ &&
$t_{ave}$ & $t_{min}$ & $t_{max}$ &&
$t_{ave}$ & $t_{min}$ & $t_{max}$ \\
\hline
10 &&     4.52 &     3.93 &     4.97 &&     0.14 &     0.13 &     0.20 &&     3.17 &     2.71 &     5.72 \\ \hline
20 &&     4.49 &     3.98 &     5.07 &&     0.18 &     0.17 &     0.22 &&     6.77 &     5.22 &    12.48 \\ \hline
30 &&     5.36 &     4.57 &     6.16 &&     0.22 &     0.20 &     0.25 &&    12.62 &    10.39 &    19.43 \\ \hline
40 &&     6.19 &     5.44 &     7.07 &&     0.28 &     0.25 &     0.32 &&    19.49 &    15.40 &    29.15 \\ \hline
50 &&     6.38 &     5.60 &     7.11 &&     0.33 &     0.29 &     0.39 &&    29.33 &    22.33 &    41.29 \\ \hline
\end{tabular}
\end{table}

\subsection{Comparing CRM models with existing models using real market data}\label{subsec:numer_stock}
\begin{figure}
\FIGURE
{\includegraphics*[width=\textwidth]{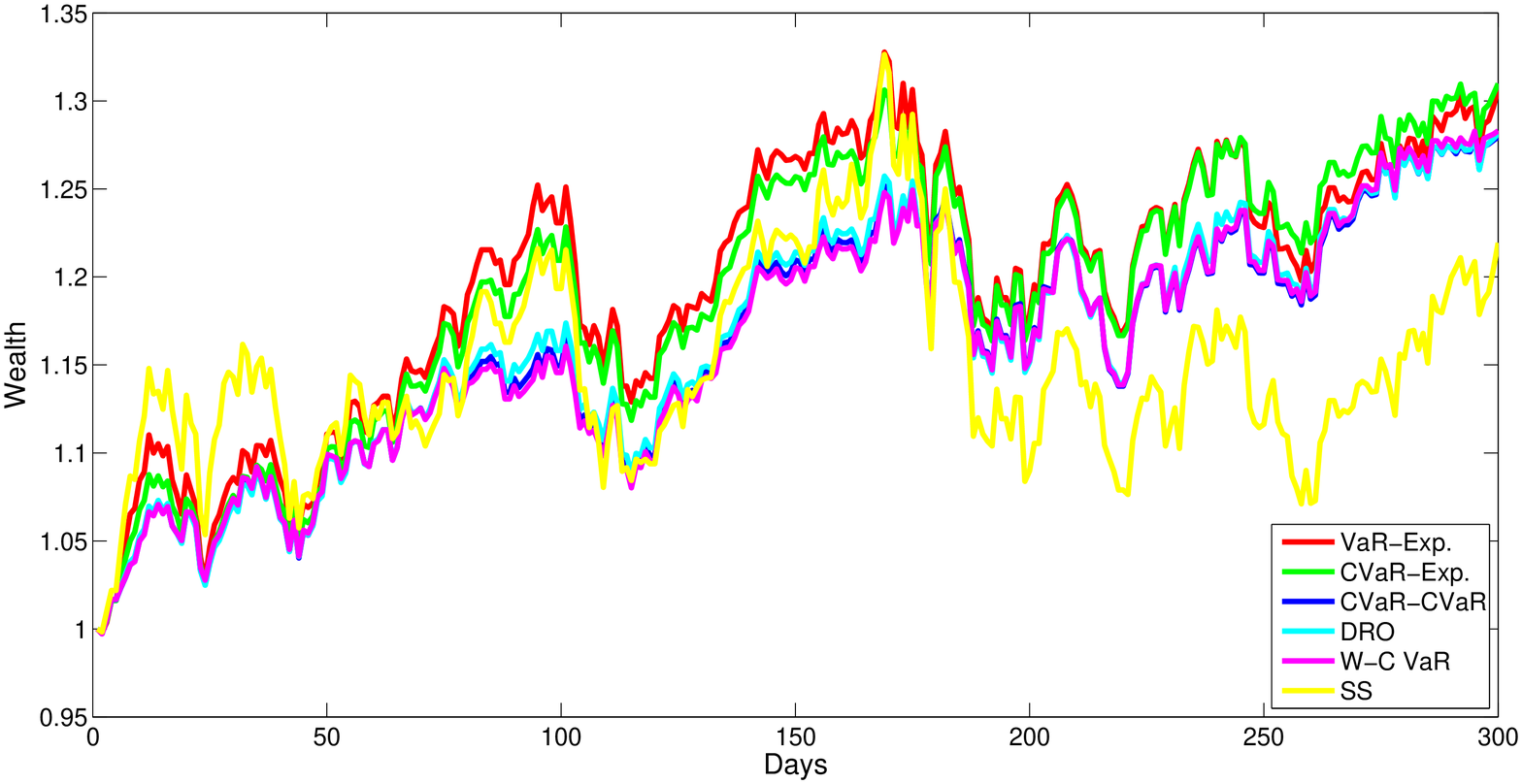}}
{Comparison of wealth accumulation of all models.\label{figure_1}}
{} 
\end{figure}

To test the performance of our new models in real world trading operations, we compare our models with existing models for portfolio selection problem. In each experiment, we randomly choose 4 stocks from all 359 stocks to build a dynamic portfolio during the period from 3/4/2010 to 4/27/2011 (300 days in total). The portfolio is recomputed everyday using the returns of the last 30 days as input data.
At each day of the experiment, only the returns of the last 30 days can be used. We set sample size as 2000 for VaR-Expectation, CVaR-Expectation and CVaR-CVaR models, and compare the results with distributionally robust model (\citealt{Delage-Ye-2010}), worst-case VaR model (\citealt{worst_case_VaR}) and single stock (SS) model. The single stock model chooses the stock that has the highest average return rate in the last 30 days as the sole stock for that day, and is used as a naive benchmark here. The average cumulative return of each day is shown in Figure \ref{figure_1}, while the average standard deviation and daily return are presented in Table \ref{experiment_table3}.
In this experiment, the volatilities of our models are significantly smaller than that of the SS approach.
The CVaR-Expectation model and VaR-Expectation model achieve better return rates than the other models, while CVaR-CVaR model performs as good as DRO model and worst-case VaR model.
Though we cannot draw general conclusions about which model is intrinsically better without more intensive tests, this experiment shows that the CRM models choose portfolios with robust performance, which is a property we desire.

\begin{table}\footnotesize \caption{Average return rate and standard deviation of all models.}
\label{experiment_table3}


\centering
\begin{tabular}{ccccccc}
\hline
&
\multicolumn{1}{c}{VaR-Exp.}&
\multicolumn{1}{c}{CVaR-Exp.}&
\multicolumn{1}{c}{CVaR-CVaR}&
\multicolumn{1}{c}{DRO}&
\multicolumn{1}{c}{W-C VaR}&
\multicolumn{1}{c}{SS}\\
\hline
ave return&
0.096 \% &      0.096 \% &      0.087 \% &      0.088 \% &      0.087 \% &      0.078 \% \\ \hline
ave std&
$ 1.49 \times 10^{-2}$ & $ 1.34 \times 10^{-2}$ & $ 1.17 \times 10^{-2}$ & $ 1.19 \times 10^{-2}$ & $ 1.17 \times 10^{-2}$ & $ 2.17 \times 10^{-2}$ \\ \hline
\end{tabular}
  \end{table}

\section{Conclusions}\label{sec:conclusion}
In this paper, we propose a unified framework for decision making under uncertainty using the composite of risk measures.
Our focus has been the case where the distribution of uncertain model parameters can be parameterized by a finite number of parameters, which includes a large family of problems.
The generality of our framework allows us to unify several existing models as well as to construct new models within the framework.
We show through theoretical proofs and numerical experiments that our new paradigms yield less conservative solutions, yet provide the same degree of probabilistic guarantee.

%
%
%




\bibliographystyle{ormsv080} 
\bibliography{references} 



\end{document}